# Objective Bayesian analysis under sequential experimentation[*]

Dongchu Sun [1] and James O. Berger [2]

*University of Missouri-Columbia and Duke University*

**Abstract:** Objective priors for sequential experiments are considered. Common priors, such as the Jeffreys prior and the reference prior, will typically depend on the stopping rule used for the sequential experiment. New expressions for reference priors are obtained in various contexts, and computational issues involving such priors are considered.

## Contents

## 1. Introduction

Bayesian analysis using objective or default priors has received considerable attention; cf. Datta and Mukerjee [17], Bernardo [6, 7] Berger and Bernardo [4], Berger [3], Ghosh, Delampady and Samanta [21], and references therein. The latter book, in particular, contains an excellent discussion of the issues and controversies involving objective priors, reflecting the many years of leadership of J. K. Ghosh in the field (along with his many coauthors). See, for example, [13, 20, 22, 23].

[*]Supported by the NSF Grants DMS-01-03265, SES-03-51523 and SES-07-20229, and the NIH Grant R01-MH071418.

[1]Department of Statistics, University of Missouri-Columbia, 146 Middlebush Hall, Columbia, MO 65211-6100, USA, e-mail: sund@missouri.edu; url: www.stat.missouri.edu/ dsun
[2]Department of Statistical Science, Duke University, Box 90251, Durham, NC 27708-0251, USA, e-mail: berger@stat.duke.edu; url: www.stat.duke.edu/ berger

A common objective prior is the Jeffreys prior [27], which is proportional to the square root of the determinant of the Fisher information matrix. The Jeffreys prior is quite useful for a single parameter model, but can be seriously deficient for multi-parameter models; this has led to preference for reference priors in multiparameter situations (cf. Berger and Bernardo [5] and Bernardo [7]).

Almost all results on objective priors have been for fixed sample size experiments. In practice, however, statistical experiments are often conducted sequentially, with a known stopping rule (cf. Siegmund [30] and Ghosh, Sen and Mukhopadhyay [24]). Bartholomew [2] and Geisser [19] introduced the notion that objective priors for a sequential experiment should depend on the expected stopping time. Ye [38] derived the reference prior for sequential experiments when the expected stopping time depends on the parameter of interest only. In this paper we generalize Ye's result in various directions, and provide some new computational tools for use with priors that depend on expected stopping times.

The paper is arranged as follows. Section 2 reviews the Fisher information matrix for sequential experiments with a known stopping rule, derives the Jeffreys/reference prior for illustrative one-parameter examples, and then provides an expression for multiparameter reference priors when the stopping rule satisfies a certain property. In Section 3, reference priors and matching priors (cf. Datta and Mukerjee [17]) are derived for Bar-Lev and Reiser's [1] two-parameter exponential family. Illustrations are given for normal distributions with several commonly used stopping times.

Computation of expected stopping times is often difficult, so that utilization of reference priors for sequential experiments is typically challenging. In Section 4, an approximation to the reference prior for sequential experiments is introduced which is exact under some circumstances, seems to be a reasonable approximation in general, and allows for much simpler computation.

## 2. Noninformative priors with a known stopping rule

### 2.1. Notation and the Jeffreys-rule prior

We assume that $X_1, X_2, \ldots,$ is an i.i.d. sequence of random variables with density $f(x \mid \boldsymbol{\theta})$ that is *regular* (Walker [35]). Here $\boldsymbol{\theta}$ is a $q \times 1$ vector of unknown parameters. Let $N$ denote a proper stopping time for the sequential experiment – see Govindarajulu [25] for a definition, which also is a source for the following well-known lemma:

**Lemma 2.1.** *Let $\mathbf{I}(\boldsymbol{\theta})$ be the Fisher information matrix based on $X_1$. Under the proper stopping time $N$, the Fisher information based on $(X_1, \ldots, X_N)$ is*

$$\mathbf{I}^* = E_{\boldsymbol{\theta}}(N)\mathbf{I}(\boldsymbol{\theta}). \tag{2.1}$$

The Jeffreys-rule prior [27] for $\boldsymbol{\theta}$ is defined as the square root of the determinant of the Fisher information matrix. In the fixed sample size case, this is $\pi_J(\boldsymbol{\theta}) \propto |\mathbf{I}(\boldsymbol{\theta})|^{1/2}$. For the sequential experiment, it follows from the above lemma that Jeffreys' prior is

$$\pi_J^*(\boldsymbol{\theta}) \propto \{E_{\boldsymbol{\theta}}(N)\}^{q/2}|\mathbf{I}(\boldsymbol{\theta})|^{1/2} \propto \{E_{\boldsymbol{\theta}}(N)\}^{q/2}\pi_J(\boldsymbol{\theta}). \tag{2.2}$$

**Example 2.1.** Let $N_r$ be a random variable with a negative binomial distribution $NB(r, p)$, where $r$ is a positive integer and $p \in (0, 1)$. Let $X_1, X_2, \ldots$ be a sequence of Bernoulli random variables with success probability $p$. $N_r$ can be viewed as a stopping time for this Bernoulli sequence as follows:

$$N_r = \inf\{n \geq 1 : X_1 + \cdots + X_n = r\}.$$

The probability of $N_r$ is

$$P(N_r = k) = \binom{k-1}{r-1} p^r (1-p)^{k-r}, \quad \text{for } k = r, r+1, \ldots.$$

An easy computation yields $E_p(N_r) = r/p$. Since the Jeffreys rule prior for a Bernoulli random variable is $\pi_J(p) \propto 1/\sqrt{p(1-p)}$, it follows from (2.2) that the Jeffreys rule prior for the negative binomial distribution is

$$\pi_J^*(p) \propto \frac{r}{p} \pi_J(p) \propto \frac{1}{p\sqrt{1-p}}.$$

This, of course, is well known from a direct computation with the negative binomial distribution, as discussed in Geisser [8] and Bernardo and Smith ([19], Example 5.14, p. 315).

We next consider an example with a continuous stopping time.

**Example 2.2.** Let $\{Z(t) : t > 0\}$ be a Brownian motion with constant drift $\theta$ and variance 1 per unit time, so $Z(t) \sim N(\theta t, t)$. Let $-\infty < a < 0 < b < \infty$, and let $T_{ab}$ denote the random stopping time

(2.3) $$T_{ab} = \inf\{t > 0 : Z(t) \leq a \text{ or } Z(t) \geq b\}.$$

It follows from Hall [26] that

$$E_\theta(T_{ab}) = \begin{cases} \frac{1}{\theta}\left[b - (b-a)\frac{e^{2b\theta}-1}{e^{2(b-a)\theta}-1}\right], & \text{if } \theta \neq 0, \\ -ab, & \text{if } \theta = 0. \end{cases}$$

Note that the constant prior is the Jeffreys prior based on stopping at a fixed time (Polson and Roberts [29]; Sivaganesan and Lingam [31] ), from which it follows that the Jeffreys or reference prior for this situation is

$$\pi(\theta) = \sqrt{E_\theta(T_{ab})}.$$

This is of additional interest because of the study in Brown [10], which showed that the commonly used estimate $Z(T)/T$, which is the posterior mean under a constant prior for $\theta$, is inadmissible under estimation with squared error loss. Brown [10] further suggested that prior distributions which behaved like $|\theta|^{-1}$ as $|\theta| \to \infty$ were optimal for this situation. The Jeffreys/reference prior has behavior $|\theta|^{-1/2}$ as $|\theta| \to \infty$, and so is not of this form, but admissibility is very dependent on the loss function used. Indeed, it can be argued that a weighted-squared error loss is appropriate for this situation, and the reference prior is likely admissible for an appropriate weight.

## 2.2. Reference priors

Reference priors depend on a grouping and ordering of the parameters; see Berger and Bernardo [4, 5]. Suppose that $\boldsymbol{\theta} = (\boldsymbol{\theta}_{(1)}, \ldots, \boldsymbol{\theta}_{(m)})$ is an $m$-ordered grouping, where the dimension of component $\boldsymbol{\theta}_{(i)}$ is $q_i$ for $i = 1, \ldots, m$. Datta and Ghosh [14] considered the special case in which the (fixed sample size) Fisher information matrix is diagonal, with the diagonal elements being products of functions of the $\boldsymbol{\theta}_{(i)}$. Our first result is a generalization of their result.

**Theorem 2.1.** *Suppose that the Fisher information matrix corresponding to a single observation $X_1$ is of the form*

$$(2.4) \qquad \mathbf{I}(\boldsymbol{\theta}) = diag\Big(\prod_{i=1}^{m} \boldsymbol{G}_{1i}(\boldsymbol{\theta}_{(i)}), \ldots, \prod_{i=1}^{m} \boldsymbol{G}_{mi}(\boldsymbol{\theta}_{(i)})\Big),$$

*where $\boldsymbol{G}_{li}$ is a $q_i \times q_i$ matrix. Assume further that the expected stopping time is of the form*

$$(2.5) \qquad E_{\boldsymbol{\theta}}(N) = \prod_{i=1}^{m} g_i(\boldsymbol{\theta}_{(i)}).$$

*Then the reference prior for $\boldsymbol{\theta}$ in the sequential experiment is*

$$(2.6) \qquad \pi_R^*(\boldsymbol{\theta}_{(1)}, \ldots, \boldsymbol{\theta}_{(m)}) \propto \prod_{i=1}^{m} [g_i(\boldsymbol{\theta}_{(i)})]^{q_i/2} \pi_R(\boldsymbol{\theta}_{(1)}, \ldots, \boldsymbol{\theta}_{(m)}),$$

*where $\pi_R(\boldsymbol{\theta}_{(1)}, \ldots, \boldsymbol{\theta}_{(m)})$ is the reference prior based on the single observation $X_1$, given by*

$$(2.7) \qquad \pi_R(\boldsymbol{\theta}_{(1)}, \ldots, \boldsymbol{\theta}_{(m)}) = \prod_{i=1}^{m} |\boldsymbol{G}_{ii}(\boldsymbol{\theta}_{(i)})|^{1/2}.$$

*Proof.* The proof is essentially identical to that in Datta [12], noting that, under (2.5), the sequential Fisher information matrix has the product structure of Datta and Ghosh [13, 14, 15]. □

This theorem can also be considered to be a generalization of Ye [38], who considered the case where $E_{\boldsymbol{\theta}}(N)$ depends only on $\boldsymbol{\theta}_{(1)}$, the parameter of interest.

Berger and Bernardo [5] suggested that one should always try to use a one-at-a-time reference prior, where each component of the grouping of parameters contains only one parameter, and much of the subsequent literature has validated this suggestion. We thus take it as given here that a one-at-a-time reference prior is the desired target. The following result is an immediate corollary of Theorem 2.1.

**Corollary 2.1.** *Suppose that the conditions of Theorem 2.1 hold. If $q_i = 1$, for $i = 1, \ldots, m = k$, then the resulting one-at-a-time reference prior for $\boldsymbol{\theta}$ in the sequential experiment is*

$$\pi_R^*(\boldsymbol{\theta}) \propto \sqrt{E_{\boldsymbol{\theta}}(N)} \pi_R(\boldsymbol{\theta}_1, \ldots, \boldsymbol{\theta}_k).$$

For later purposes, we also note another corollary of Theorem 2.1, which applies if the dimension of each component of the grouping of parameters has dimension 2.

**Corollary 2.2.** *Suppose that the conditions of Theorem 2.1 hold. If all $q_j = 2$, then the reference prior for $\boldsymbol{\theta}$ in the sequential experiment is*

$$\pi_R^*(\boldsymbol{\theta}) \propto E_{\boldsymbol{\theta}}(N) \pi_R(\boldsymbol{\theta}_{(1)}, \ldots, \boldsymbol{\theta}_{(m)}).$$

## 3. A two-parameter exponential family

### 3.1. The model and reference priors

Bar-Lev and Reiser [1] considered the following density function of the generic two-parameter exponential family:

$$(3.1) \quad f(x \mid \theta_1, \theta_2) = a(x) \exp\{\theta_1 U_1(x) - \theta_1 G'_2(\theta_2) U_2(x) - \psi(\theta_1, \theta_2)\},$$

where $\theta_1 < 0$, $\theta_2 = E\{U_2(X) \mid (\theta_1, \theta_2)\}$, $G_i(\cdot), (i = 1, 2)$ are infinitely differentiable functions satisfying $G''_i > 0$, and

$$\psi(\theta_1, \theta_2) = -\theta_1\{\theta_2 G'_2(\theta_2) - G_2(\theta_2)\} + G_1(\theta_2).$$

This is a large class of distributions, which includes, for suitable choices of $G_1$, $G_2$, $U_1$ and $U_2$, many popular statistical models such as the normal, inverse normal, gamma, and inverse gamma. Table 1, reproduced from Sun [32], indicates how each distribution arises.

Let $X_1, X_2, \ldots$ be a sequence of random variables from (3.1). The Fisher information per observation is

$$\mathbf{I}(\theta_1, \theta_2) = \begin{pmatrix} G''_1(\theta_1) & 0 \\ 0 & -\theta_1 G''_2(\theta_2) \end{pmatrix}.$$

The two parameters $\theta_1$ and $\theta_2$ are orthogonal in the sense of Cox and Reid [11]. Thus the Jeffreys prior for a single observation is

$$(3.2) \quad \pi_J(\theta_1, \theta_2) \propto \sqrt{|\theta_1|} \sqrt{G''_1(\theta_1) G''_2(\theta_2)}.$$

When either $\theta_1$ or $\theta_2$ is the parameter of interest, it is shown in Sun and Ye [33] that the one-at-a-time reference priors are

$$(3.3) \quad \pi_R(\theta_1, \theta_2) = \sqrt{G''_1(\theta_1) G''_2(\theta_2)}.$$

The parameter $\theta_2$ is the expectation of $U_2(X_1)$. Bose and Boukai [9] considered inference about $\theta_2$ in sequential experimentation with the following stopping time:

$$(3.4) \quad N_a = \inf\left\{n \geq m_0 : Y_n < nG'_1\left(-\frac{a^2}{n^2}\right)\right\}, \quad a \geq 0,$$

where $Y_n = n^{-1} \sum_{i=1}^n U_1(X_i) - G_2\{n^{-1} \sum_{i=1}^n U_2(X_i)\}$ and $m_0 \geq 2$ is an initial sample size. From Theorem 2 of Bose and Boukai [9], we have

$$(3.5) \quad \lim_{a \to \infty} \frac{N_a}{a} = \frac{1}{\sqrt{|\theta_1|}} \quad a.s.$$

$$(3.6) \quad \lim_{a \to \infty} E_{\boldsymbol{\theta}}\left(\frac{N_a}{a}\right) = \frac{1}{\sqrt{|\theta_1|}}.$$

TABLE 1
*Special cases of Bar-Lev and Reiser's [1] two parameter exponential family, where* $h(\theta_1) = -\theta_1 + \theta_1 \log(-\theta_1) + \log(\Gamma(-\theta_1))$

| | $G_1(\theta_1)$ | $G_2(\theta_2)$ | $U_1(x)$ | $U_2(x)$ | $\theta_1$ | $\theta_2$ |
|---|---|---|---|---|---|---|
| $N(\mu, \sigma^2)$ | $-\frac{1}{2}\log(-2\theta_1)$ | $\theta_2^2$ | $x^2$ | $x$ | $-1/(2\sigma^2)$ | $\mu$ |
| Inverse Gaussian | $-\frac{1}{2}\log(-2\theta_1)$ | $1/\theta_2$ | $1/x$ | $x$ | $-\alpha/2$ | $\sqrt{\alpha/\mu}$ |
| Gamma | $h(\theta_1)$ | $-\log\theta_2$ | $-\log x$ | $x$ | $-\alpha$ | $\mu$ |
| Inverse Gamma | $h(\theta_1)$ | $-\log\theta_2$ | $\log x$ | $1/x$ | $-\alpha$ | $\mu$ |

Bar-Lev and Reiser [1] showed that the distribution of $Y_n$ does not depend on the parameter $\theta_2$. So condition (2.5) satisfies when either $\theta_1$ or $\theta_2$ is the parameter of interest. The following result is immediate from Theorem 2.1 or Corollary 2.1.

**Fact 3.1.** (a) The Jeffreys prior for $(\theta_1, \theta_2)$ in model (3.1) with the stopping time (3.4) and when $a$ is large is approximately

$$\pi_J^*(\theta_1, \theta_2) \propto \sqrt{G_1''(\theta_1) G_2''(\theta_2)}. \tag{3.7}$$

(b) The one-at-a-time reference prior for $(\theta_1, \theta_2)$ in model (3.1), when either $\theta_1$ or $\theta_2$ is the parameter of interest, the stopping time (3.4) is used, and $a$ is large enough, is approximately

$$\pi_R^*(\theta_1, \theta_2) \propto \frac{1}{|\theta_1|^{1/4}} \sqrt{G_1''(\theta_1) G_2''(\theta_2)}. \tag{3.8}$$

**Example 3.1.** Suppose $X_1, X_2, \ldots,$ are a sequence of $N(\mu, \sigma^2)$ random variables. Then $\theta_1 = -1/2\sigma^2, \theta_2 = \mu$, $G_1'(\theta_1) = -1/2\theta_1$, and $Y_n = \sum_{i=1}^n (X_i - \overline{X}_n)^2$. The stopping rule (3.4) becomes

$$N_a = \inf\Big\{n \geq m_0 : n^{-1} \sum_{i=1}^n (X_i - \overline{X}_n)^2 < n^2/(2a^2)\Big\}.$$

So the priors (3.2), (3.3), (3.7), and (3.8) are, respectively,

$$\pi_J(\mu, \sigma^2) \propto \frac{1}{(\sigma^2)^{3/2}}, \ \pi_R(\mu, \sigma^2) \propto \frac{1}{\sigma^2}, \ \pi_J^*(\mu, \sigma^2) \propto \frac{1}{\sigma^2}, \ \pi_R^*(\mu, \sigma^2) \propto \frac{1}{(\sigma^2)^{3/4}}$$

or equivalently,

$$\pi_J(\mu, \sigma) \propto \frac{1}{\sigma^2}, \ \pi_R(\mu, \sigma) \propto \frac{1}{\sigma}, \ \pi_J^*(\mu, \sigma) \propto \frac{1}{\sigma}, \ \pi_R^*(\mu, \sigma) \propto \frac{1}{\sqrt{\sigma}}.$$

### 3.2. Probability matching priors for a sequential experiment

Asymptotic frequentist coverage is an often-used criterion to compare objective priors; see Welch and Peers [36], Peers [28], Tibshirani [34], Datta and Ghosh [13], Datta, Ghosh and Mukerjee [16], and Datta and Mukerjee [17] for discussion and references. The most common approach is to find a "matching prior," i.e., a prior which results in posterior one-sided credible intervals that are also accurate as frequentist confidence intervals. Another type of matching prior, considered by Sun and Ye [33], is a prior such that the confidence interval based on the signed squared root transformation of the log-likelihood ratio is also a Bayesian credible interval. Almost all of the literature considers the fixed sample case for i.i.d. observations; exceptions are Ye [38] and Sun [32].

For sequential experiments involving the Bar-Lev and Reiser [1] two-parameter exponential family, let $l_n(\theta_1, \theta_2)$ be the log-likelihood function of $(\theta_1, \theta_2)$, given $\boldsymbol{X}_n = (X_1, \ldots, X_n)$, and let $(\hat\theta_{n1}, \hat\theta_{n2})$ be the maximum likelihood estimator of $(\theta_1, \theta_2)$. Write

$$Y_n = n^{-1} \sum_{i=1}^n U_1(X_i) - G_2\Big\{n^{-1} \sum_{i=1}^n U_2(X_i)\Big\}.$$

Then, on $\{Y_n \in G'_1(\Theta_1)\} \cap \{n^{-1} \sum_{i=1}^n U_2(X_i) \in \Theta_2\}$, $\hat{\theta}_{n1}$ is the solution of $Y_n = G'_1(\hat{\theta}_{n1})$, and $\hat{\theta}_{n2} = n^{-1} \sum_{i=1}^n U_2(X_i)$. Define

$$I_1(\omega_1, \theta_1) = G_1(\theta_1) - G_1(\omega_1) - G'_1(\omega_1)(\theta_1 - \omega_1), \quad \omega_1, \theta_1 \in \Theta_1,$$
$$I_2(\omega_2, \theta_2) = G_2(\omega_2) - G_2(\theta_2) - G'_2(\theta_2)(\omega_2 - \theta_2), \quad \omega_2, \theta_2 \in \Theta_2.$$

From the convexity of $G_1$ and $G_2$, these two functions are nonnegative. From Sun [32], the log-likelihood ratio is $l_n(\hat{\theta}_{n1}, \hat{\theta}_{n2}) - l_n(\theta_1, \theta_2) = (Z_{n1}^2 + Z_{n2}^2)/2$, where

$$\begin{pmatrix} Z_{n1} \\ Z_{n2} \end{pmatrix} = \begin{pmatrix} \{2nI_1(\hat{\theta}_{n1}, \theta_1)\}^{1/2} \, sgn(\theta_1 - \hat{\theta}_{n1}) \\ \{-2n\theta_1 I_2(\hat{\theta}_{n2}, \theta_2)\}^{1/2} \, sgn(\theta_2 - \hat{\theta}_{n2}) \end{pmatrix}$$

is a generalized signed square root of the log-likelihood ratio.

Let $P_{(\theta_1,\theta_2)}$ denote probability over $X_1, X_2, \ldots$, given $(\theta_1, \theta_2)$, and, for a fixed prior $\pi(\theta_1, \theta_2)$, let $P^\pi(\cdot \mid \boldsymbol{X}_n)$ denote posterior probability given $\boldsymbol{X}_n$. Suppose we are considering a stopping time, $N_a$, such that $N_a \to \infty$ almost surely as $a \to \infty$. An asymptotic frequentist matching prior in this sequential setting is a prior $\pi$ such that

$$(3.9) \qquad P^\pi(Z_{N_a,1} \leq c_1, Z_{N_a,2} \leq c_2 \mid \boldsymbol{X}_{N_a})$$
$$= P_{(\theta_1,\theta_2)}(Z_{N_a,1} \leq c_1, Z_{N_a,2} \leq c_2) + O(a^{-1}),$$

for all $c_1$ and $c_2$ in $P_{(\theta_1,\theta_2)}$−probability.

Suppose now that the stopping rule satisfies

$$(3.10) \qquad \frac{N_a}{a} \to \tau(\boldsymbol{\theta}), \quad \text{in } L_1.$$

From Sun [32], the unique prior satisfying (3.9), and hence the unique asymptotic matching prior, is

$$(3.11) \qquad \pi_m^*(\theta_1, \theta_2) \propto \sqrt{\tau(\boldsymbol{\theta}) G''_1(\theta_1) G''_2(\theta_2)}.$$

As an immediate example, for the stopping time defined in (3.4), property (3.6) establishes that (3.10) holds; hence the reference prior given in (3.8) is also the asymptotic matching prior, a very desirable situation.

**Example 3.1** (continued). In deriving the sequential likelihood ratio test to see if $(\mu, \sigma^2) = (\mu_0, \sigma_0^2)$, Woodroofe [37] considered the following stopping rule,

$$(3.12) \qquad N_a = \min\left(b_2 a, \ \inf\left\{n \geq b_1 a : \sum_{i=1}^n X_i^2 - n - n\log(\hat{\sigma}_n^2) > 2a\right\}\right),$$

where $0 < b_1 < b_2 < \infty$ are two prespecified numbers, $\hat{\sigma}_n^2 = n^{-1} \sum_{i=1}^n (X_i - \overline{X}_n)^2$, and $\overline{X}_n = n^{-1} \sum_{i=1}^n X_i$. Theorem 8.3 of Woodroofe [37] implies that

$$\frac{a}{N_a} \to \begin{cases} b_2, & \text{if } \rho^2(\boldsymbol{\theta}) < 1/b_2, \\ \rho^2(\boldsymbol{\theta}), & \text{if } 1/b_2 < \rho^2(\boldsymbol{\theta}) < 1/b_1, \\ b_1, & \text{if } \rho^2(\boldsymbol{\theta}) > 1/b_1, \end{cases}$$

in $P_{(\theta_1,\theta_2)}$−probability, as $a \to \infty$, where

$$(3.13) \qquad \rho^2(\boldsymbol{\theta}) = G_1(\theta_1) - G_1(-0.5) - G'_1(-0.5)(\theta_1 + 0.5) - \theta_1 \theta_2^2$$
$$= \{(\mu^2 + 1)/\sigma^2 + \log(\sigma^2) - 1\}/2.$$

Thus (3.11) gives an asymptotic matching prior for this situation. Note, however, that the expected stopping time is not of the form (2.5), so that we cannot assert that this prior is also a one-at-a-time reference prior.

## 4. Computation

If $E_{\boldsymbol{\theta}}[N]$ is available in closed form, as in the examples in this paper, computation with any of the sequential priors can be done using common MCMC techniques. Hence we only consider here the case in which $E_{\boldsymbol{\theta}}[N]$ can only be computed numerically.

### 4.1. Brute force computation

All the Jeffreys, reference, and matching priors that have been discussed for a sequential experiment are of the form $\Psi(E_{\boldsymbol{\theta}}[N])\pi_F(\boldsymbol{\theta})$, where $\Psi$ is some operator and $\pi_F$ is the corresponding prior for the fixed sample size experiment. The posterior distribution corresponding to this prior is

$$(4.1) \qquad \pi^*(\boldsymbol{\theta} \mid \boldsymbol{X}_N) \;\propto\; \Psi(E_{\boldsymbol{\theta}}[N])\pi_F(\boldsymbol{\theta}) \prod_{i=1}^{N} f(X_i \mid \boldsymbol{\theta}),$$

where $\boldsymbol{X}_N = (X_1, \ldots, X_N)$ is the data.

The brute force method for simulating from this posterior distribution is the following Metropolis algorithm:

*Step 1.* Sample a proposed $\boldsymbol{\theta}'$, from the fixed sample size posterior density of $\boldsymbol{\theta}$, which is proportional to $\pi_F(\boldsymbol{\theta}) \prod_{i=1}^{N} f(X_i \mid \boldsymbol{\theta})$.

*Step 2.* Numerically estimate $E_{\boldsymbol{\theta}'}[N]$. For instance, one could repeatedly sample $N$ from its distribution given $\boldsymbol{\theta}'$, by simply repeatedly simulating the sequential experiment for the given $\boldsymbol{\theta}'$, observing the $N$ that results from each simulation, and averaging to obtain the estimate $\widehat{E_{\boldsymbol{\theta}'}[N]}$.

*Step 3.* Perform a Metropolis step: sample $u \sim$ uniform $(0,1)$ and, with $\boldsymbol{\theta}$ denoting the previous value the parameter, accept $\boldsymbol{\theta}'$ if

$$u \leq \min\left\{1, \frac{\Psi(\widehat{E_{\boldsymbol{\theta}}[N]})}{\Psi(\widehat{E_{\boldsymbol{\theta}'}[N]})}\right\},$$

and set $\boldsymbol{\theta}'$ equal to the previous $\boldsymbol{\theta}$ otherwise.

If one cannot directly draw from the posterior in Step 1, one could instead using any MCMC scheme, e.g. Gibbs sampling or Metropolis–Hastings. If doing so, however, be sure to iterate Step 1 many times before moving on to Step 2. This is because Step 2 is typically extremely expensive, as it may involve thousands of simulations of the entire experiment simply to compute one Metropolis acceptance probability. In situations where one dependent step is much more expensive than others, it pays to iterate first on the others.

### *4.2. The two-dimensional case*

If using the Jeffreys prior in a two-dimensional problem or the reference prior in the situation of Corollary 2.2, the posterior distribution is of the form

$$(4.2) \qquad \pi^*(\boldsymbol{\theta} \mid \boldsymbol{X}_N) \quad \propto \quad E_{\boldsymbol{\theta}}[N] \, \pi_F(\boldsymbol{\theta}) \prod_{i=1}^{N} f(X_i \mid \boldsymbol{\theta}) \,.$$

This allows a remarkable simplification in the computation, by introducing $N$ as a latent variable.

To avoid confusion, we will label the latent variable as $\tilde{N}$; it is a variable with the same distribution as $N$, but is independent of $N$. Write the density of $\tilde{N}$ given $\boldsymbol{\theta}$ as $p(\tilde{N} \mid \boldsymbol{\theta})$. Then the joint density of $(\tilde{N}, \boldsymbol{\theta})$, given the data $\boldsymbol{X}_N = (X_1, \ldots, X_N)$, is proportional to

$$(4.3) \qquad p(\tilde{N} \mid \boldsymbol{\theta}) \tilde{N} \, \pi_F(\boldsymbol{\theta}) \prod_{i=1}^{N} f(X_i \mid \boldsymbol{\theta}) \,.$$

Sampling $(\tilde{N}, \boldsymbol{\theta})$ from this distribution will result in $\boldsymbol{\theta}$ from (4.2), as can easily be seen by marginalizing over $\tilde{N}$ in (4.3).

Here is a Metropolis algorithm for sampling from (4.3).

*Step 1.* Sample a proposed $\boldsymbol{\theta}'$, from the fixed sample size posterior density of $\boldsymbol{\theta}$, which is proportional to $\pi_F(\boldsymbol{\theta}) \prod_{i=1}^{N} f(X_i \mid \boldsymbol{\theta})$.

*Step 2.* Sample a proposed $\tilde{N}'$ from $p(\tilde{N} \mid \boldsymbol{\theta}')$. This can always be done by simply simulating the sequential experiment once, given $\boldsymbol{\theta}'$.

*Step 3.* Perform a Metropolis step: sample $u \sim \text{uniform}(0,1)$ and, letting $(\tilde{N}, \boldsymbol{\theta})$ denote the previous value the parameter, accept $(\tilde{N}', \boldsymbol{\theta}')$ if

$$u \leq \min\left\{1, \frac{\tilde{N}}{\tilde{N}'}\right\},$$

and set $(\tilde{N}', \boldsymbol{\theta}')$ equal to the previous $(\tilde{N}, \boldsymbol{\theta})$ otherwise. (Note that, if $\tilde{N}' < \tilde{N}$, one would always accept the candidate.)

The reason that this is a vastly more efficient algorithm than the brute force algorithm is that one need only simulate a single draw of $\tilde{N}'$ in Step 2, whereas thousands of draws would be needed in Step 2 of the brute force algorithm to compute $\widehat{E_{\boldsymbol{\theta}'}[N]}$. Again, of course, Step 1 could be replaced by any convenient dependent MCMC scheme. Whether one then needs to iterate Step 1 before moving on to Step 2 will be context dependent.

### *4.3. Modified reference priors*

The most desirable prior is the one-at-a-time reference prior given in Corollary 2.1, resulting in the posterior distribution

$$(4.4) \qquad \pi^*(\boldsymbol{\theta} \mid \boldsymbol{X}_N) \quad \propto \quad \sqrt{E_{\boldsymbol{\theta}}[N]} \, \pi_R(\boldsymbol{\theta}) \prod_{i=1}^{N} f(X_i \mid \boldsymbol{\theta}) \,.$$

Unfortunately, the latent variable trick is not available for sampling from this distribution.

Interestingly, however, it is frequently the case that

(4.5) $$\sqrt{E_{\boldsymbol{\theta}}[N]} \approx E_{\boldsymbol{\theta}}[\sqrt{N}].$$

When this is the case, the latent variable trick can be applied, and the algorithm from Section 4.2 can be utilized by simply replacing $\tilde{N}/\tilde{N}'$ in the Metropolis step with $\sqrt{\tilde{N}/\tilde{N}'}$.

In the remainder of the section, we discuss the reason that the approximation (4.5) often holds. The first is that the sampling distribution of $N$ may be rather concentrated in a region of large $N$, in which case the approximation is clearly good.

**Example (Bar-Lev and Reiser [1]) (continued).** For the stopping time $N_a$ defined in (3.4), it follows from (3.5) and (3.6) that

$$\lim_{a\to\infty} E_{\boldsymbol{\theta}}(\sqrt{N_a/a}) = 1/|\theta_1|^{1/4}.$$

We then have

$$\lim_{a\to\infty} E_{\boldsymbol{\theta}}\sqrt{\frac{N_a}{a}} \propto \lim_{a\to\infty} \sqrt{E_{\boldsymbol{\theta}}\left(\frac{N_a}{a}\right)}.$$

**Example 2.1 (continued).** Let $N_r$ have the negative binomial distribution $NB(r, p)$. Note that $E_p(N_r) = r/p$ and $Var_p(N_r) = rp/(1-p)^2$. As $r \to \infty$, we have

$$\sqrt{N_r/r} \to 1/\sqrt{p} \text{ in probability}$$

and

$$E_p(\sqrt{N_r/r}) \to \sqrt{E_p(N_r/r)} \equiv 1/\sqrt{p}.$$

To see the difference between $E_p(\sqrt{N_r/r})$ and $\sqrt{E_p(N_r/r)}$ for moderate $r$, they are plotted, as a function of $p$, in Figure 1 for $r = 1$ and $r = 9$. For $r = 9$, the curves are essentially indistinguishable; even for the minimal $r = 1$ they are quite close.

It is also interesting to look at the posterior distributions for this example. In Figure 2, we plot the posterior densities of $p$ for three priors

$$\begin{aligned}\pi_J(p) &\propto 1/\sqrt{p(1-p)},\\ \pi_R^*(p) &\propto 1/\sqrt{p},\\ \pi_M(p) &\propto E_p(\sqrt{N_r^*/r}).\end{aligned}$$

Here $\pi_M(p)$ is an approximate prior. For even the very small $r = 2$, the posterior densities under the two priors $\pi_R^*$ and $\pi_M$ are quite close, yet substantially different from that under $\pi_J$. For a moderate $r = 10$, the posterior densities under $\pi_R^*$ and $\pi_M$ are almost identical. Note that the posterior densities of $p$ under $\pi_J$ and $\pi_R^*$ are Beta $(r, N_r - r + 0.5)$ and Beta $(r = 0.5, N_r - r + 0.5)$, respectively. The posterior densities of $p$ under $\pi_M$ were computed using 5000 Metropolis samples.

As a final indication of the similarity of the true and approximate reference priors in this example, and of the value of using the sequential reference priors, we compare the frequentist coverage probabilities that result from their use in obtaining confidence intervals for $p$. Table 2 considers the frequentist coverage of one-sided 5% and 95% Bayesian credible regions, based on the fixed sample size Jeffreys' prior $\pi_J$, the sequential Jeffreys/reference prior $\pi_R^*$ and the approximate prior $\pi_M$ for various combination of $r$ and $p$. The fixed sample size Jeffreys' prior performs worse then the other two, indicating the value of using the sequential versions, while the reference prior and the approximate prior are almost equally good.

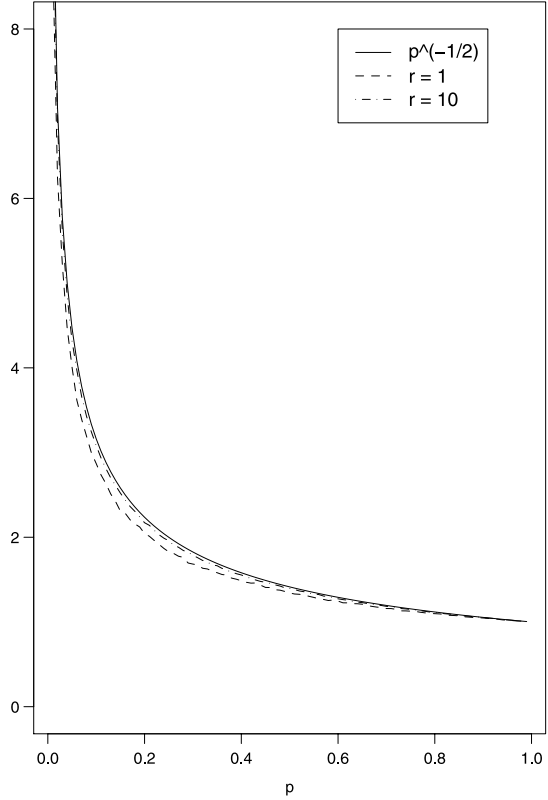

FIG 1. *Negative binomial example: comparison of $\sqrt{E_p(N_r/r)}$ and $E_p(\sqrt{N_r/r})$ for $r = 1$ and $r = 9$.*

TABLE 2
*Coverage Probability of one-sided 5% (95%) Bayesian credible sets for the negative binomial Example 2.1, under the three priors $\pi_J(p) = 1/\sqrt{p(1-p)}$, $\pi_R^*(p) = 1/(p\sqrt{1-p})$, and*
$$\pi_M(p) = E_p(\sqrt{N_r^*/r})$$

| $r$ | $p$ | $\pi_J$ | $\pi_R^*$ | $\pi_M$ |
|---|---|---|---|---|
| 2  | 0.1 | 0.1142(0.9738) | 0.0516(0.9511) | 0.0487(0.9509) |
| 2  | 0.5 | 0.0002(0.9652) | 0.0010(0.9381) | 0.0008(0.9455) |
| 2  | 0.9 | 0.0001(0.9724) | 0.0003(0.9700) | 0.0000(0.9729) |
| 8  | 0.1 | 0.0751(0.9642) | 0.0474(0.9498) | 0.0465(0.9534) |
| 8  | 0.5 | 0.0552(0.9688) | 0.0522(0.9536) | 0.0568(0.9517) |
| 8  | 0.9 | 0.0000(0.9307) | 0.0001(0.9310) | 0.0002(0.9339) |
| 30 | 0.1 | 0.0617(0.9571) | 0.0508(0.9497) | 0.0516(0.9523) |
| 30 | 0.5 | 0.0556(0.9594) | 0.0512(0.9495) | 0.0525(0.9503) |
| 30 | 0.9 | 0.0426(0.9369) | 0.0438(0.9410) | 0.0442(0.9368) |

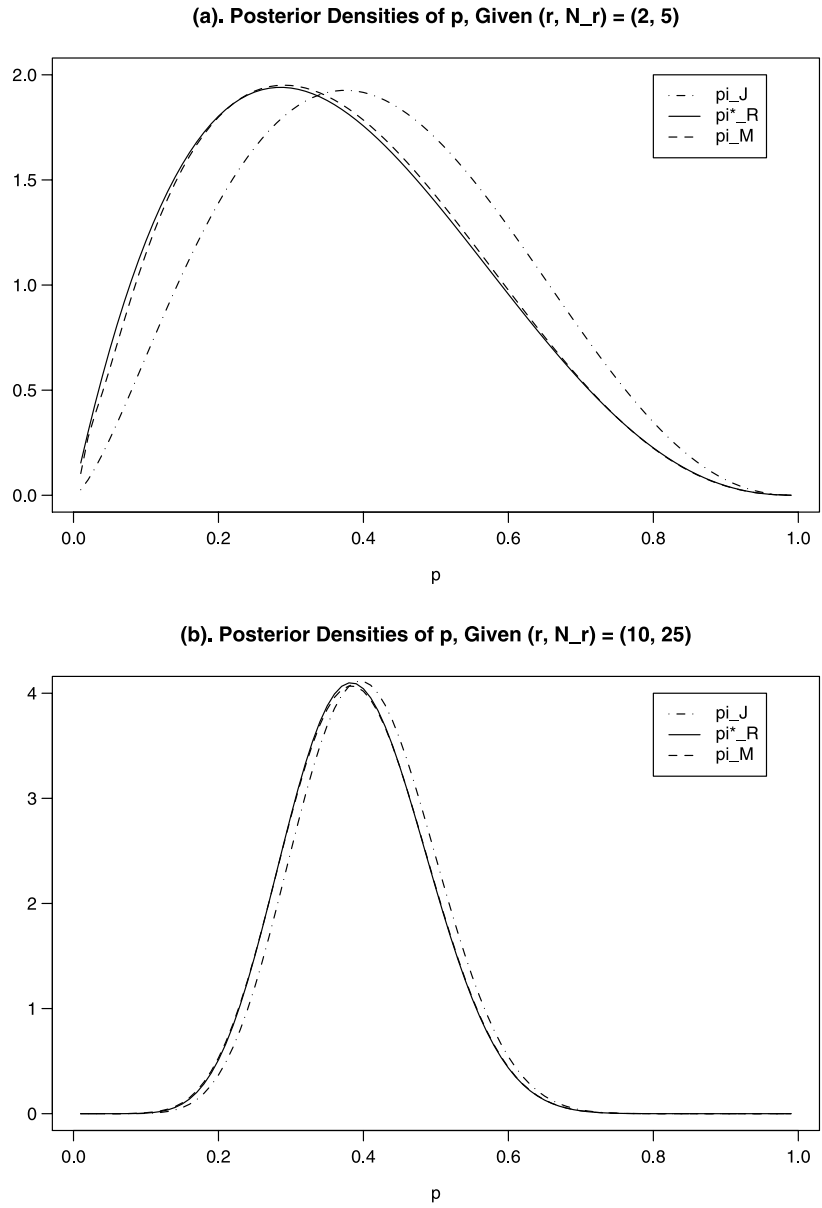

FIG 2. *Posterior densities of $p$ based on the priors $\pi_J(p) = 1/\sqrt{p(1-p)}$, $\pi_R^*(p) = 1/(p\sqrt{1-p})$, and $\pi_M(p) = E_p(\sqrt{N_r^*/r})$ for $r = 1, 10$; (a) $(r, N_r) = (2, 5)$; (b) $(r, N_r) = (10, 25)$.*

**Acknowledgments.** The authors gratefully acknowledge the comments and suggestions of a referee.